\documentclass[11pt]{article}

\usepackage{mathrsfs}
\usepackage{latexsym}
\usepackage{amsmath,amsthm}
\usepackage{amsfonts}
\usepackage{url}

\usepackage{amssymb}

\newcommand{\h}{\eta}

\usepackage{bbm}
\newcommand{\chr}{\boldsymbol{\mathbbm{1}}} 
\newcommand{\pred}[1]{\chr_{\left\{ #1 \right\}}}
\newcommand{\prs}{\vec{P}}
\newcommand{\pr}[1]{\prs\!\tlprn{#1}}
\newcommand{\mexp}{\vec{E}}

\newcommand{\E}{\mexp}
\renewcommand{\P}{\prs}
\newcommand{\Pn}{\P\!_n}
\newcommand{\Pnn}[1]{\P\!_{#1}}

\newcommand{\Dn}{\Delta_n}
\newcommand{\Gn}{\Gamma_n}
\newcommand{\trn}{^{\!\mathsf{T}}} 
\newcommand{\inv}{^{-1}} 

\newcommand{\TV}[1]{\nrm{#1}_{\textrm{{\tiny \textup{TV}}}}}
\newcommand{\ham}{\textrm{{\tiny \textup{Ham}}}}
\newcommand{\Lip}[1]{\nrm{#1}_{\textrm{{\tiny \textup{Lip}}}}}
\newcommand{\Lipp}[2]{\nrm{#1}_{\textrm{{\tiny \textup{Lip}}},#2}}
\newcommand{\Lips}{\operatorname{Lip}}

\newcommand{\sgn}{\operatorname{sgn}}
\newcommand{\diam}{\operatorname{diam}}
\newcommand{\essup}{\mathop{\operatorname{ess\,sup}}}

\newcommand{\one}{\text{i}}
\newcommand{\two}{\text{ii}}
\newcommand{\thr}{\text{iii}}

\newcommand{\subplus}{_{
+
}}

\newcommand{\pl}[1]{\paren{#1}\subplus}

\usepackage{geometry}
\makeatletter
\theoremstyle{plain}
\newtheorem{thm}{Theorem}[section]
\theoremstyle{plain}    
\newtheorem*{thm*}{Theorem}
 \theoremstyle{remark}
 \newtheorem{rem}[thm]{Remark}
 \theoremstyle{plain}    
 \newtheorem{lem}[thm]{Lemma} 
\makeatother
\newcommand{\bethn}{\begin{thm}}
\newcommand{\enthn}{\end{thm}}
\renewcommand{\beth}{\begin{thm*}}
\newcommand{\enth}{\end{thm*}}
\newcommand{\bepf}{\begin{proof}}
\newcommand{\enpf}{\end{proof}}
\newcommand{\belen}{\begin{lem}}
\newcommand{\enlen}{\end{lem}}
\newcommand{\ben}{\begin{enumerate}}
\newcommand{\een}{\end{enumerate}}
\newcommand{\bit}{\begin{itemize}}
\newcommand{\eit}{\end{itemize}}

\renewcommand{\vec}[1]{\bs{\mathrm{#1}}}

\newcommand{\phinorm}[1]{\nrm{#1}_{\Phi}}
\newcommand{\psinorm}[1]{\nrm{#1}_{\Psi}}
\newcommand{\basicspace}{{\cal S}}
\newcommand{\X}{\basicspace}
\newcommand{\supr}[1]{^{(#1)}}

\newcommand{\seq}[3]{#1_{#2}\ldots#1_{#3}}
\newcommand{\sseq}[3]{#1_{#2}^{#3}}  
\newcommand{\cd}{\,}
\newcommand{\ccat}[1]{[#1]}
\newcommand{\sd}{}
\newcommand{\scat}[1]{#1}

\newcommand{\dsabs}[1]{\bigl| #1 \bigr|}

\newcommand{\nrm}[1]{\left\Vert #1 \right\Vert}

\newcommand{\tsnrm}[1]{\Vert #1 \Vert}
\newcommand{\iprod}[2]{\left\langle #1 , #2 \right\rangle}

\newcommand{\f}{\varphi}
\renewcommand{\k}{\kappa}
\newcommand{\el}{\ell}

\newcommand{\psin}[1]{\Psi_{#1}}

\newcommand{\calA}{\mathcal{A}}
\newcommand{\calB}{\mathcal{B}}

\newcommand{\calF}{\mathcal{F}}

\newcommand{\calL}{\mathcal{L}}

\newcommand{\calX}{\mathcal{X}}

\newcommand{\tha}{\theta}
\newcommand{\RR}{\mathbb{R}}
\newcommand{\NN}{\mathbb{N}}
\newcommand{\ZZ}{\mathbb{Z}}
\newcommand{\QQ}{\mathbb{Q}}
\newcommand{\beq}{\begin{eqnarray*}}
\newcommand{\eeq}{\end{eqnarray*}}
\newcommand{\beqn}{\begin{eqnarray}}
\newcommand{\eeqn}{\end{eqnarray}}
\newcommand{\paren}[1]{\left( #1 \right)}
\newcommand{\sqprn}[1]{\left[ #1 \right]}
\newcommand{\tlprn}[1]{\left\{ #1 \right\}}
\newcommand{\set}[1]{\tlprn{#1}}
\newcommand{\abs}[1]{\left| #1 \right|}
\newcommand{\ceil}[1]{\ensuremath{\left\lceil#1\right\rceil}}
\newcommand{\gn}{\, | \,}

\newcommand{\ts}{\textstyle}
\newcommand{\bs}{\boldsymbol}

\def\blk{~\texttt{<LK>}~}
\def\elk{~\texttt{</LK>}~}
\newcommand{\lnote}[1]{\blk #1 \elk}
\newcommand{\hide}[1]{}
\newcommand{\oo}[1]{\frac{1}{#1}}

\def\eps{\varepsilon}
\newcommand{\defeq}{\doteq}
\newcommand{\abscont}{\ll}

\title{
Metric and Mixing
Sufficient Conditions \\
for
Concentration
of
Measure 
}
\author{Leonid Kontorovich\\
School of Computer Science\\ 
Carnegie Mellon University\\ 
Pittsburgh, PA 15213\\
USA \\
\url{lkontor@cs.cmu.edu}
}
\begin{document}
\maketitle
\abstract{
We derive sufficient conditions for a family $(\X^n,\rho_n,\Pn)$ of metric probability
spaces to 
have the measure concentration property.
Specifically, if the sequence
$\{\Pn\}$ of probability measures satisfies a strong mixing condition
(which we call $\eta$-mixing) and the sequence of metrics $\{\rho_n\}$
is
what we call
$\Psi$-dominated, we show that $(\X^n,\rho_n,\Pn)$ is a 
normal L\'evy family. We establish these properties for some metric
probability spaces, including the possibly novel $\X=[0,1]$,
$\rho_n=\nrm{\cdot}_1$ case.
}


\section{Introduction}
\subsection{Background}
\label{sec:intro}
The study of measure concentration in general metric spaces was
initiated 
in the 1970's
by Vitali Milman,
who in turn drew inspiration from
Paul L\'evy's work
(see~\cite{talagrand95} for a brief historical exposition). Since
then, various deep insights have been gained into the concentration of
measure phenomenon~\cite{ledoux01}.

The words ``measure'' and ``concentration'' suggest an interplay of
analytic and 
geometric
aspects. Indeed, there are two essential
ingredients in proving a concentration result: the 
random variable
must be continuous in a strong (Lipschitz) sense, and the
random process 
must be mixing in some
strong sense. The simple examples 
we give
in \S\ref{sec:conc-martingale} 
illustrate how,
in general, the failure of either of these conditions to hold
can prevent 
a random variable
from being concentrated.

A common way of summarizing the
phenomenon is to say that in a high-dimensional space, almost all of
the probability is concentrated 
around
any set whose measure is at least
${\ts\oo2}$. Another way is to say that any
``sufficiently continuous'' function is tightly concentrated about
its mean. 
To state this more formally (but still somewhat imprecisely), let 
$(X_i)_{1\leq i\leq n}$,
$X_i\in\X$,
be the random process defined on
the probability space $(\X^n,\calF,\P)$,
and $f:\X^n\to\RR$ be a function satisfying some Lipschitz condition
(and possibly others, such as convexity). A concentration of measure
result 
(for our purposes)
is an inequality of the form
\beqn
\label{eq:basicdev}
\pr{\abs{f(X)-\E f(X)}>t} &\leq& c\exp(-K t^2)
\eeqn
where $c>0$ is a small constant (typically, $c=2$) and $K>0$ is
some quantitative indicator of the strong mixing properties of $X$. It
is crucial that neither $c$ nor $K$ depend on $f$.\footnote{See
\cite{ledoux01} for a much more general notion of concentration.}

A few celebrated milestones that naturally fall into the 
paradigm 
of (\ref{eq:basicdev})
include
L\'evy's original isoperimetric inequality
on the 
sphere (see the notes and references in
\cite{LedouxTal91}), McDiarmid's bounded differences
inequality~\cite{mcd89}, and Marton's generalization of~\cite{mcd89}
for contracting Markov chains~\cite{marton96}.
(Talagrand's no-less celebrated series of results~\cite{talagrand95}
does not easily lend itself to such a compact description.)

Building on the work of Azuma~\cite{azuma} and Hoeffding~\cite{hoeff},
McDiarmid showed that 
if $f:\X^n\to\RR$
has $\Lip{f}\leq1$
under the normalized Hamming metric $\bar d_\ham$
and $\P$ is a product measure on $\X^n$, we have
\beqn
\label{eq:mcd}
\pr{\abs{f-\E f}>t} &\leq& 2\exp(-2nt^2)
\eeqn
(he actually proved this for the more general class of weighted
Hamming metrics).
Using coupling and information-theoretic inequalities, Marton showed
that if the conditions on $f:\X^n\to\RR$ are as above and $\P$ is a
contracting Markov measure on $\X^n$ with Doeblin coefficient $\tha<1$,
\beqn
\label{eq:marton}
\pr{\abs{f-M_f}>t} &\leq &
2\exp\sqprn{
-2n\paren{
t(1-\tha)
-\sqrt{\frac{\log2}{2n}}}^2
},
\eeqn
where $M_f$ is a $\P$-median of $f$. Since product measures are
degenerate cases of Markov measures (with $\tha=0$), Marton's result
is a powerful generalization of (\ref{eq:mcd}).

Two natural directions for extending results of type (\ref{eq:mcd}) are
to 
derive
such inequalities for various measures (processes) and
metrics. Talagrand's paper
\cite{talagrand95}
is a tour de force in proving concentration for various (not
necessarily metric) notions of distance, but it deals exclusively with
product measures.
Since 
the publication of
Marton's
concentration inequality 
in 1996 
(to our knowledge, 
the first of its kind
for a
nonproduct, non-Haar measure),
several authors proceeded to generalize her information-theoretic
approach~\cite{dembo97,dembo-zeitouni},
and offer alternative approaches based on the entropy method
\cite{ledoux96,samson00}
or martingale techniques~\cite{kontram06}.
Talagrand in~\cite{talagrand95} discusses strengths and weaknesses of
the martingale
method, observing that
``while in principle the martingale method has a wider range
of applications, in many situations the [isoperimetric] inequalities
[are] more powerful.'' Bearing out his first point, Kontorovich and
Ramanan~\cite{kontram06} used martingales to derive a general strong
mixing condition for concentration (in the $\bar d_\ham$ metric),
applying it to weakly contracting Markov chains. Following up, 
Kontorovich
extended the technique to hidden Markov
\cite{kont06-hmm}
 and Markov tree
\cite{kont06-tree}
measures.

Although a detailed survey of measure concentration literature 
is not our intent
here, we remark that
many of the results mentioned above may be
described as 
working
to
extend inequalities of
type (\ref{eq:basicdev}) to 
wider classes of measures
and metrics by imposing
different strong mixing and Lipschitz continuity conditions. Already
in~\cite{marton96}, Marton gives a (rather stringent) mixing condition
sufficient for concentration. Later, Marton~\cite{marton03,marton04}
and Samson~\cite{samson00} prove concentration for general classes of
processes in terms of various mixing coefficients; Samson applies this
to Markov chains and $\phi$-mixing processes while Marton's
application concerns lattice random fields.

In this paper, we build upon the results in~\cite{kontram06} and give
general metric and mixing conditions 
that
ensure 
the
concentration
of
measure. 
We make use of a fundamental mixing coefficient, 
which
has appeared (under various guises) in Marton's and Samson's work, to
define the notion of $\eta$-mixing for a random process. We also
define a condition on the metric space, which we call $\Psi$-dominance.
\hide{
them may be cast in this paradigm
(some recent papers explicitly refer to {\it mixing} in the title,
\cite{marton03,samson00}). In this paper, we introduce 
\lnote{define, introduced in K+R...}
a strong mixing
condition called $\eta$-mixing, similar in spirit to the quantities
defined in \cite{marton03,samson00}. Our other condition is
on the metric space, namely that of being
$\Psi$-dominated. 
}
Our main result,
Theorem~\ref{thm:main},
states that if the family of metric probability spaces
$(\X^n,\rho_n,\P)_{n\geq 1}$ is such that $\P$ is $\eta$-mixing and 
$(\X^n,\rho_n)_{n\geq1}$ is $\Psi$-dominated, then
$(\X^n,\rho_n,\P)$ is a normal L\'evy family, and therefore exhibits
measure concentration. We also give examples of metric probability
spaces satisfying these conditions.

\subsection{Paper outline}
\label{sec:outline}
This paper is organized as follows. In
\S\ref{sec:notconv}, we fix some notation used throughout the paper and
dispose of some measure-theoretic issues. We 
review
L\'evy families and concentration functions, and their connection to
deviation inequalities in \S\ref{sec:levy}. 
In \S\ref{sec:conc-martingale} we introduce the method of bounded
martingale differences as our technique for proving measure
concentration. We define the two key notions of this paper,
{\em $\eta$-mixing} and
{\em $\Psi$-dominance} in 
\S\ref{sec:hmix} and
\S\ref{sec:psidom}, respectively.
Our main concentration result
for $\eta$-mixing processes with $\Psi$-dominated metrics is proved in
\S\ref{sec:main}.
In \S\ref{sec:pdex} we give examples of some natural $\Psi$-dominated
metrics,
and conclude the paper with a summary and brief discussion in \S\ref{sec:disc}.
Finally, the Appendix takes a bit of a scenic detour,
examining the
two norms defined in this paper and
the strength of the topologies they induce.

\section{Notation and technicalities}
\label{sec:notconv}
Random variables are capitalized ($X$), specified sequences (vectors)
are written in lowercase ($x\in\X^n$), the shorthand
$\sseq{X}{i}{j}\defeq
(X_i,\ldots,X_j)
$ is used for all sequences, and
brackets denote sequence concatenation:
$\ccat{\sseq{x}{i}{j}\cd\sseq{x}{j+1}{k}}=\sseq{x}{i}{k}$.
Often,
for readability, we abbreviate
$\ccat{y\cd w}$ as
$yw$.

We use the indicator variable
$\pred{\cdot}$ 
to assign 0-1 truth values 
to the predicate in 
$\set{\cdot}$.
The sign function is defined by
$\sgn(z)=\pred{z>0}-\pred{z<0}$. The ramp function
is defined by $\pl{z}=z\pred{z>0}$.

We will follow Talagrand's time-honored tradition of dispensing with
measure-theoretic technicalities, since the (well-understood) problems
they raise would distract us from the big picture.
Only in the Appendix do these issues become interesting and relevant,
and are handled there with rigor.

\hide{
We will use generic measure space notation $(\calX,\calF,\mu)$ when
talking about abstract metric probability spaces and L\'evy families;
when using the language of deviation inequalities we will write 
$(\X^n,\calF,\Pn)$ to emphasize the correspondence between the measure
$\Pn$ and the random process $(X_i)_{1\leq i\leq n}$.
}
In any metric probability space
$(\calX,\rho,\P)$,
it is 
understood that 
$\P$ is a measure on
the Borel 
$\sigma$-algebra generated from the topology induced by $\rho$.
We 
will often
abuse notation slightly by
suppressing the dependence on the dimensionality $n$ in 
the measures
$\Pn$.
In such cases, 
we are
implicitly assuming that the 
probability
measures are {\em consistent} in the sense that
for each Borel set $A\subset \X^{n-1}$, we have
\beq
\label{eq:consist}
\Pnn{n-1}(A)
= \int_{A\times\X} d\Pn(\sseq{x}{1}{n}).
\eeq
The probability 
$\P$
and expectation $\E$ operators
are defined 
with respect 
the measure space specified in context.
To any probability space $(\X^n,\calF,\P)$, we associate the canonical
random process $X=\sseq{X}{1}{n}$, $X_i\in\X$, satisfying
\beq
\pr{X\in A} &=& \P(A)
\eeq
for any $A\in\calF$.

If $\mu$ is a positive Borel measure on $(\calX,\calF)$ and $\tau$ is a
signed measure on $(\calX,\calF)$, 
we define
the {\em total variation} of $\tau$ 
by
\beqn
\label{eq:tv}
2\TV{\tau} &=& \sup\sum_{i=1}^\infty\abs{\tau(E_i)},
\eeqn
where the supremum is over all the countable partitions ${E_i}$ of
$\calX$ (this 
quantity
is necessarily finite, by Theorem 6.4 of~\cite{rudin}).\footnote{
Note the factor of $2$ in (\ref{eq:tv}), which typically does not
appear in analysis texts but is standard in probability theory, when
$\tau$ is the difference of two probability measures.}
It is a consequence of the Lebesgue-Radon-Nikod\'ym theorem
(\cite{rudin}, Theorem 6.12) that if
$\tau\abscont\mu$ with density $h$, we have
\beq
2\TV{\tau} &=& \int_{\calX} |h|d\mu.
\eeq
Additionally, if $\tau$ is
{\em balanced}, meaning that $\tau(\calX)=0$, we have
\beqn
\label{eq:taubal}
\TV{\tau} &=& \int_{\calX} \pl{h} d\mu;
\eeqn
this follows 
the Hahn decomposition 
(\cite{rudin}, Theorem 6.14).

If $(\calX,\calF,\mu)$ is a (positive) measure space, we write $L_p(\calX,\mu)$ for the
usual space of $\mu$-measurable functions $f:\calX\to\RR$, whose $L_p$
norm
\beq
\nrm{f}_{
L_p(\calX,\mu)
} &=& \paren{\int_\calX |f|^p d\mu}^{1/p}
\eeq
is finite. We will write 
$\nrm{\cdot}_{L_p(\calX,\mu)}$
as
$\nrm{\cdot}_{L_p(\mu)}$
or just
$\nrm{\cdot}_{L_p}$ if there is no ambiguity;
when $\mu$ is the counting measure on a discrete space, we write this
as $\nrm{\cdot}_p$.

Likewise, the $L_\infty$ norm,
$\nrm{f}_{L_\infty}=\essup\abs{f}$
 is defined via the essential supremum:
\beq
\essup_{x\in\calX}f(x)
&=&
\inf\{a\in[-\infty,\infty]:\mu\{f(x)>a\}=0\}.
\eeq

The Hamming metric on a product space $\X^n$ is the sum of the
discrete metrics on 
$\X$:
\beq
d_\ham(x,y) &=& \sum_{i=1}^n \pred{x_i\neq y_i}
\eeq
for $x,y\in\X^n$. Sometimes we will work with the normalized Hamming
metric: $\bar d_\ham = \oo n d_\ham$.

\section{
L\'evy families
and concentration
in metric spaces
}
\label{sec:levy}
A natural language for discussing measure concentration in general
metric spaces is that of L\'evy families.
This 
definition
is taken, with minor variations,
from Chapter 6
of~\cite{milman-schechtman}. Let $(\calX,\rho,\P)$ be a 
Borel 
probability
space whose topology is induced by the metric $\rho$.
Whenever we write $A\subset \calX$, it is implicit that $A$ is a Borel
subset of $\calX$. For $t>0$, define the $t$-{\em fattening} of
$A\subset\calX$:
$$ A_t = \{ x\in\calX: \rho(x,A)\leq t\}.$$
The {\em concentration function} 
$\alpha(\cdot)=\alpha_{
\calX,\rho,\P
}(\cdot)
$
is defined by:
\beq
\alpha(t) &=& 1-\inf\{\P(A_t): A\subset\calX, ~\P(A)\geq{\ts\oo2}\}.
\eeq
Let $(\calX_n,\rho_n,\Pn)_{n\geq 1}$
be a family of metric
probability spaces with $\diam_{\rho_n}(\calX_n)<\infty$,
where
\beqn
\diam_{\rho_n} (\calX_n) &\defeq& \sup_{x,y\in\calX_n} \rho_n(x,y).
\eeqn
 This
family\hide{
\footnote{
Whenever we omit the subscript $n$ from $\mu_n$ or $\Pn$, we are
implicitly assuming that the measures are {\em consistent} in the sense that
for $A\subset \X^{n-1}$, we have
$$ \mu_{n-1}(A) = \int_{A\times\calX} d\mu_n(\sseq{x}{1}{n}).$$
}}
is called a 
{\em normal L\'evy} family if there are constants $c_1,c_2>0$ such
that
\beq
\alpha
_{\calX_n,\rho_n,\Pn}
(t
) &\leq& c_1\exp(-c_2 nt^2)
\eeq
for each $t>0$ and $n\geq1$.
\hide{
\begin{rem}
\label{rem:mu_n}
Whenever we omit the subscript $n$ from $\mu_n$ or $\Pn$, we are
implicitly assuming that the measures are {\em consistent} in the sense that
for $A\subset \X^{n-1}$, we have
$$ \mu_{n-1}(A) = \int_{A\times\X} d\mu_n(\sseq{x}{1}{n}).$$
\end{rem}
}

The condition of being a normal L\'evy family implies strong
concentration of a Lipschitz $f:\calX_n\to\RR$ about its median (and
mean); this connection is explored in-depth in~\cite{ledoux01}. In
particular, if $(\calX,\rho,\P)$ is a metric probability space and
$f:\calX\to\RR$ is measurable, define 
its {\em modulus of continuity} by
\beqn
\label{eq:modcont}
\omega_f(\delta) &=& \sup\{\abs{f(x)-f(y)} : \rho(x,y)<\delta \}.
\eeqn
A number $M_f\in\RR$ is called a {\em median} of $f$ if
\beq
\pr{f\leq M_f}\geq{\ts\oo2}
\qquad\text{and}\qquad
\pr{f\geq M_f}\geq{\ts\oo2}
\eeq
(a median need not be unique). These definitions immediately imply
the deviation inequality
\cite{ledoux01}(1.9)
\beq
\pr{\abs{f-M_f}>\omega_f(\delta)} &\leq& 2\alpha
_{\calX,\rho,\P}
(\delta),
\eeq
which in turn yields~\cite{ledoux01}(1.13)
\beqn
\label{eq:dev-ledoux}
\pr{\abs{f-M_f}>t} &\leq& 2\alpha
_{\calX,\rho,\P}
(t/\Lip{f}),
\eeqn
where the Lipschitz constant $\Lip{f}$ is the smallest 
constant
$C$ for
which $\omega_f(\delta)\leq C\delta$, for all $\delta>0$.
In particular, (\ref{eq:dev-ledoux}) lets us take $\Lip{f}=1$ without
loss of generality, which we shall do below.
Theorem 1.8 in~\cite{ledoux01} lets us convert concentration about a
median to concentration about any constant:
\beth
Let $f$ be a measurable function on a probability space
$(\calX,\calA,\P)$.
Assume that for some $a\in\RR$ and a non-negative function
$\alpha$ on $\RR_+$ such that $\lim_{r\to\infty}\alpha(r)=0$,
\beq
\pr{{\abs{f-a}\geq r}} &\leq& \alpha(r)
\eeq
for all $r>0$. Then
\beq
\pr{\abs{f-M_f}\geq    r+r_0} &\leq& \alpha(r),\quad r>0,
\eeq
where $M_f$ is a 
$\P$-median of $f$ and where $r_0>0$ is
such that $\alpha(r_0)<{\ts\oo2}$.
If moreover 
$\bar\alpha=\int_0^\infty \alpha(r)dr<\infty$ then
$f$ is integrable, $\abs{a-\E f}\leq\bar\alpha$, and for every
$r>0$, 
\beq
\pr{\abs{f-\E f}\geq    r+\bar\alpha} &\leq& \alpha(r).
\eeq
\enth
Thus, 
for a normal L\'evy family, deviation inequalities for the mean and
median are equivalent up to the constants $c_1,c_2$. Theorem 1.7 in
\cite{ledoux01} is a converse to (\ref{eq:dev-ledoux}), showing that if
Lipschitz functions on a metric probability space 
$(\calX,\rho,\P)$
are tightly
concentrated about their means, this implies a rapid decay
of
$\alpha_{\calX,\rho,\P}(\cdot)$.

\hide{
\lnote{Many measure concentration results amount to establishing that
  a certain 
family of metric
probability spaces is a normal L\'evy family. Discussion in Ledoux
  p. 55; Talagrand[check!!], Marton, K+R fall into
  this category; check Saurav, Dembo, Bobkov. Not Samson, since
  requires convexity}
}

\section{Concentration via martingale differences}
\label{sec:conc-martingale}

\subsection{Background}
\label{sec:mgale-bkg}
Let 
$(\X^n,\calF,\P)$
be a probability space, where
$\calF$ is the usual Borel $\sigma$-algebra generated by the finite
dimensional cylinders.
On this space
define the random process
$(X_i)_{1\leq i\leq n}$, 
$X_i\in\X$.
Let $\calF_i$ be the $\sigma$-algebra generated by 
$(\seq{X}{1}{i})$, 
which induces the filtration
\beqn
\label{eq:filtr}
\{\emptyset,\X^n\} = \calF_0 \subset\calF_1 \subset\ldots
\subset\calF_n = \calF.
\eeqn
For $i=1,\ldots,n$ and 
$f\in L_1(\X^n,\P)$,
define the martingale difference
\beqn
\label{eq:mdifdef}
V_i &=& \E[f\gn \calF_i] - \E[f\gn \calF_{i-1}].
\eeqn
It is a classical result,\footnote{
See~\cite{ledoux01} for a modern presentation and a short proof of (\ref{eq:azuma}).
}
going back to Azuma~\cite{azuma}, that
\beqn
\label{eq:azuma}
\pr{\abs{f-\E f}>t} &\leq& 2\exp(-t^2/2D^2)
\eeqn
where $D^2\geq\sum_{i=1}^n\nrm{V_i}_\infty^2$
(the meaning of $\nrm{V_i}_\infty$ will be made explicit later). 
Thus, if we are able to
uniformly bound the martingale difference,
\beq
\max_{1\leq i\leq n} \nrm{V_i}_\infty &\leq& H_n,
\eeq
we obtain the concentration inequality
\beqn
\label{eq:genconc}
\pr{\abs{f-\E f}>t} &\leq& 2\exp
\paren{-\frac{t^2}{2 n H_n^2}}
.
\eeqn
Our ability to derive results of the type in (\ref{eq:genconc}) will
in general
depend on the continuity properties of $f$ and the mixing properties
of the process $X$.

Let us give two simple examples to build up some intuition. Let $\P$
be the uniform probability measure on $\{0,1\}^n$ and
$(X_i)_{1\leq i\leq n}$ be the associated (independent)
process. Though 
different notions of mixing exist~\cite{bradley},
$X$ trivially satisfies them all, being an i.i.d. process. Define
$f:\{0,1\}^n\to[0,1]$ 
by
$$ f(x) = 
x_1 \oplus x_2 \oplus\ldots\oplus x_n
,
$$
where $\oplus$ is addition mod 2. Since
$\pr{f(X)=0}=\pr{f(X)=1}={\ts\oo2}$, 
$f$ is certainly not concentrated
about its mean (or any other constant).
Though $X$ is as well-behaved as can be,
$f$
is ill-behaved in the sense that flipping any single
input bit causes the output to fluctuate by $1$.\footnote{
Without making far-reaching claims, we comment on a possible
connection between the
oscillatory behavior of $f$ and
the notorious
difficulty of learning noisy parity functions~\cite{feld}. By
contrast,\hide{
This oscillatory behavior of $f$ might 
be related to
the notorious
difficulty of learning noisy parity functions~\cite{feld}, while
}
the
problem of learning conjunctions and disjunctions under noise
has been solved some time ago~\cite{kearns}.}
\hide{
\begin{rem}
This example shows that McDiarmid's bound 
(\ref{eq:mcd}), recast here as
\beqn
\label{eq:mcdmod}
\pr{\abs{f-\E f}>t} &\leq& 
2\exp
\paren{-\frac{2t^2}{n\Lip{f}^2}}
\eeqn
 is not sharp
in the sense that in certain instances it vacuously bounds deviation
probabilities by a quantity that exceeds 1 for all $n$. In this
example, we have $\Lip{f}=1$ and so
(\ref{eq:mcdmod}) bounds
$\pr{\abs{f-\E f}\geq\oo2}$ 
by
$2\exp(-1/2n)>1$, for all $n$.
\end{rem}}

For the second example, take 
$f:\{0,1\}^n\to[0,1]$ 
to be
$$ f(x) = 
\oo n
\sum_{i=1}^n x_i.$$
If $(X_i)_{1\leq i\leq n}$ is the i.i.d. process from the previous
example, it is easy to show that the martingale difference in
(\ref{eq:mdifdef}) is bounded by 
$1/n$,
and so by (\ref{eq:genconc}), $f$
is concentrated about its mean. What if we relax the independence
condition? The simplest 
kind of dependence in a random process is Markovian.
Consider the homogeneous Markov process: $\pr{X_1=0}=\pr{X_1=1}=\oo2$
and $X_{i+1}=X_i$ with probability 1. This process trivially fails to
satisfy any
(reasonable) definition of mixing~\cite{bradley}.
Our well-behaved $f$ is no longer concentrated, since
we again have
$\pr{f(X)=0}=\pr{f(X)=1}={\ts\oo2}$.

\newcommand{\yi} {\sseq{y}{1}{i-1}}
\newcommand{\yii}{\sseq{y}{1}{i}}
\newcommand{\si} {\sseq{x}{1}{i-1}}
\newcommand{\sii}{\sseq{x}{1}{i}}
\newcommand{\Si} {\sseq{X}{1}{i-1}}
\newcommand{\Sii}{\sseq{X}{1}{i}}
\newcommand{\xin}{\sseq{x}{i+1}{n}}
\newcommand{\Xin}{\sseq{X}{i+1}{n}}
The two examples 
above
show that if we are to
have any hope of obtaining inequalities such as (\ref{eq:genconc}), we
will need conditions of continuity and mixing on $f$ and $X$,
respectively. Much of the discussion in 
the remainder of
this section builds upon the
treatment in~\cite{kontram06} for discrete spaces.

\subsection{Simple bound on the martingale difference}
\label{sec:simplemartbd}
Let $(\X^n,\calF,\P)$
be a probability space
and
$(X_i)_{1\leq i\leq n}$
its associated random process; define the 
filtration $\set{\calF_i}$ as in (\ref{eq:filtr}).
At this point, we make the additional assumption that
$d\P(x)=p(x)d\mu^n(x)$ for some 
positive Borel
product
measure 
$\mu^n=\mu\otimes\mu\otimes\ldots\otimes\mu$
on $(\X^n,\calF)$,
which we refer to as the {\em carrying measure}.
In the cases of interest, $\X$ will
be either countable or a compact subset of $\RR$, and correspondingly,
$\mu$ will be the counting or Lebesgue measure.
Similarly, the conditional probability
$\P(\cdot\gn\calF_i)\abscont\mu^{n-i}$, with density
$p(\cdot\gn \sseq{X}{1}{i}=\sseq{y}{1}{i})$.\hide{
For readability, we will 
sometimes
write $dx$ in the integrals instead
of the more accurate $d\mu(x)$;
also,
h}
Here and below
$p(\sseq{x}{j}{n}\gn \sseq{y}{1}{i})$ 
will occasionally be used
in place of
$p(\sseq{x}{j}{n}\gn \sseq{X}{1}{i}=\sseq{y}{1}{i})$; 
no
ambiguity should
arise.

For $f\in L_1(\X^n,\P)$,
$1\leq i\leq n$
and $\sseq{y}{1}{i}\in\X^i$, define
\beqn
\label{eq:videf}
V_i(f;\sseq{y}{1}{i}) &=&
\E[f(X)\gn\sseq{X}{1}{i}=\sseq{y}{1}{i}] -
\E[f(X)\gn\sseq{X}{1}{i-1}=\sseq{y}{1}{i-1}];
\eeqn
this is just the martingale difference. 
A slightly more tractable
quantity turns out to be
\beqn
\label{eq:vhat}
\hat V_i(f;\yi,w_i,w_i') &=&
\E[f(X)\gn\Sii = \scat{\yi \sd w_i}]-
\E[f(X)\gn\Sii = \scat{\yi \sd w_i'}],
\eeqn
where
$w_i,w_i'\in\X$.
These two quantities have a simple relationship, which may be stated
symbolically as
$\tsnrm{V_i(f;\cdot)}_{L_\infty(\P)} \leq 
\tsnrm{\hat V_i(f;\cdot)}_{L_\infty(\P)}$ and is proved in the following
lemma.
\belen
\label{lem:vhat}
Suppose $f\in L_1(\X^n,\P)$ and
and $\yii\in\X^{i}$. Then for any $\eps>0$
there are 
$w_i,w_i'\in\X$ such that
\beqn
\dsabs{V_i(f;\yii)} &\leq& \dsabs{\hat V_i(f;\yi,w_i,w'_i)}+\eps.
\eeqn
\enlen
\bepf
Let 
$$ a= 
\E[f(X)\gn\Sii = \yii]=
\int_{\X^{n-i}} 
p(\xin\gn\yii)f(\scat{\yii\sd\xin})d\mu^{n-i}(\xin);
$$
then
\beq
V_i(f;\yii) &=&
 a-
\int_{\X^{n-i+1}}
p(\sseq{x}{i}{n}\gn\yi)f(\scat{\yi\sd\sseq{x}{i}{n}})
d\mu^{n-i+1}(\sseq{x}{i}{n})\\
&=&
 a-
\int_{\X}
p(z\gn\yi) \paren{
\int_{\X^{n-i}}
p(\xin\gn\scat{\yi\sd z})
f(\scat{\yi \sd z \sd \xin})d\mu^{n-i}(\xin)}
d\mu(z)
\eeq
where the last step invokes Fubini's theorem.
We use 
the simple fact that for integrable $g,h\geq0$,
\beq
\inf h(z)\int g(z)dz \;\leq\; 
\int g(z)     h(z)dz \;\leq\; 
\sup h(z)\int g(z)dz,
\eeq
together with
$\int_\X p(z\gn\yi) d\mu(z)=1$,
to deduce, for any $\eps>0$, the 
existence of a $w_i'\in\X$ such that
\beq
\abs{V_i(f;\yii)}
&\leq&
\abs{a-
\int_{\X^{n-i}}
p(\xin\gn \scat{\yi\sd w_i'})
f(\scat{\yi\cd w_i'\sd \xin})d\mu^{n-i}(\xin)
}+\eps
\eeq
for some $w_i'\in\X$.
Taking $w_i=y_i$, this proves the claim. 
\enpf

\subsection{Martingale difference as a linear functional}
The next step is to notice that 
$\hat V_i(\cdot;\yi,w_i,w_i')$, as a functional on $L_1(\X^n,\P)$, is
linear; in fact, it is given by
\beqn
\label{eq:linfunc}
\hat V_i(f;\yi,w_i,w_i') \;=\;
\int_{\X^n} f(x) \hat g(x)d\mu^n(x)
\;\defeq\;
\iprod{f}{\hat g},
\eeqn
where
\beqn
\hat g(x) &=&
\pred{\sii= \scat{\yi\sd w_i}}
p(\xin\gn \scat{\yi\sd w_i})
-
\pred{\sii =\scat{\yi\sd w_i'}}
p(\xin\gn \scat{\yi\sd w_i'}).
\eeqn
The plan is to bound $\iprod{f}{\hat g}$ using continuity properties
of $f$ and mixing properties of $X$, which will immediately lead to a
result of type (\ref{eq:genconc}) via Lemma~\ref{lem:vhat}.

\section{$\eta$-mixing}
\label{sec:hmix}
\subsection{Definition}
\hide{
Let $(\X^n,\calF,\P)$
and
$(X_i)_{1\leq i\leq n}$
be as defined in \S\ref{sec:conc-martingale}.
and assume, as in
\S\ref{sec:simplemartbd}, 
$\P\abscont\mu$ with density 
$d\P(x)=p(x)d\mu(x)$
for some Borel 
measure
$\mu$ on $\X^n$.}
Let $(\X^n,\calF,\P)$
be a probability space
and
$(X_i)_{1\leq i\leq n}$
its associated random process.
In this section, we define
a notion of mixing particularly
suitable to our needs.
For $1\leq i<j\leq n$
and $x\in\X^i$,
let
$$\calL(\sseq{X}{j}{n}\gn \sseq{X}{1}{i}=x)
$$ be the law 
(distribution)
of $\sseq{X}{j}{n}$ conditioned on 
$\sseq{X}{1}{i}=x$. 
For $y\in\X^{i-1}$ and 
$w,w'\in\X$, 
define
\beqn
\label{eq:hdef}
\h_{ij}(y,w,w') &=&
\TV{
\calL(\sseq{X}{j}{n}\gn \sseq{X}{1}{i}=\scat{y\sd w})-
\calL(\sseq{X}{j}{n}\gn \sseq{X}{1}{i}=\scat{y\sd w'})
},
\eeqn
where $\TV{\cdot}$ is the total variation norm
(see \S\ref{sec:notconv}),
and
\beq
\bar\h_{ij} &=&
\essup_{y\in\X^{i-1},w,w'\in\X}
\h_{ij}(y,w,w'),
\eeq
where the essential supremum is taken with respect to 
the measure $\P$ on $\X^i$.
Recall 
that if 
$(U,\cal U,\P)$,
is a probability space
and
$f:U\to\RR^+$ 
is measurable,
$
\essup_{x\in U}f(x)$
is the smallest $a\in[0,\infty]$
for which
$f\leq a$ holds $\P$-almost surely.
\hide{
we have
$\essup_{x\in U}f(x)
=
\inf\{a\in[0,\infty]:\mu\{f(x)>a\}=0\}$.
}

Let $\Dn$ be the upper-triangular $n\times n$ matrix defined by
$(\Dn)_{ii}=1$ and
\beqn
\label{eq:Ddef}
(\Dn)_{ij} = \bar\h_{ij}.
\eeqn
for $1\leq i<j\leq n$. Recall that the $\ell_\infty$
operator norm is given by
\beqn
\tsnrm{\Dn}_\infty
&=&
\label{eq:infnorm}
\max_{1\leq i< n}(
1 + \bar\h_{i,i+1} + \ldots +\bar\h_{i,n}
).
\eeqn

A probability measure $\P$ on $(\X^n,\calF)$ defines the function
$H_{\P}:\NN\to\RR$ by
\beqn
\label{eq:HPn}
H_{\P}(n) &=& 
\nrm{\Dn}_\infty
;
\eeqn
we say that
the process $X$ 
(measure $\P$)
is
{\em $\eta$-mixing} if
\beqn
\label{eq:HP}
\sup_{n\to\infty} H_{\P}(n) &\defeq& \bar H_{\P} < \infty.
\eeqn

As a trivial observation, note that if the variables $(X_i)$ are
mutually independent, we have
$(\Dn)_{ij}=\pred{i=j}$ and
$\nrm{\Dn}_\infty=1$.

\subsection{Connection to $\phi$-mixing}
Samson~\cite{samson00},
using techniques quite different from those here,
showed that
if $\X=[0,1]$,
and $f:[0,1]^n\to\RR$ is
convex 
with $\Lip{f}\leq1$ (in the $\ell_2$ 
metric), 
then
\beqn
\label{eq:samson}
\pr{\abs{f(X)-\mexp f(X)}>t} &\leq& 2\exp\paren{-\frac{t^2}{2\nrm{\Gn}_2^2}}
\eeqn
where $\nrm{\Gn}_2$ is the $\el_2$ operator norm of the
matrix\footnote{
Samson used the stronger $\sup$ as opposed to $\essup$ in his analogue
of $\bar\h_{ij}$; we shall largely ignore this distinction in our analysis.
}
\beqn
\label{eq:Gndef}
(\Gn)_{ij}=\sqrt{
(\Dn)_{ij}
}.
\eeqn
Following Bradley~\cite{bradley}, 
for the random process $(X_i)_{i\in\ZZ}$ on
$(\X^\ZZ,\calF,\P)$,
we define the $\phi$-mixing
coefficient
\beqn
\phi(k) &=& \sup_{j\in\ZZ} \phi(\sseq{\calF}{-\infty}{j},\sseq{\calF}{j+k}{\infty}),
\eeqn
where $\sseq{\calF}{i}{j}\subset\calF$ is the $\sigma$-algebra generated by 
the
$\sseq{X}{i}{j}$,
and for
the $\sigma$-algebras $\calA,\calB\subset\calF$,
$\phi(\calA,\calB)$ is defined by
\beqn
\phi(\calA,\calB) &=& \sup\{\abs{\P(B\gn A)-P(B)}:
A\in\calA,~B\in\calB,~\P(A)>0\}.
\eeqn
Samson observes that 
\beqn
\bar\h_{ij}&\leq& 2\phi_{j-i},
\eeqn
 which follows from
\beq
\TV{
\calL(\sseq{X}{j}{n}\gn \sseq{X}{1}{i}=\scat{\sseq{y}{1}{i-1}\sd w})-
\calL(\sseq{X}{j}{n}\gn \sseq{X}{1}{i}=\scat{\sseq{y}{1}{i-1}\sd w'})}
&\leq&
\TV{
\calL(\sseq{X}{j}{n}\gn \sseq{X}{1}{i}=\scat{\sseq{y}{1}{i-1}\sd w})-
\calL(\sseq{X}{j}{n})}\\
&+&
\TV{
\calL(\sseq{X}{j}{n}\gn \sseq{X}{1}{i}=\scat{\sseq{y}{1}{i-1}\sd w'})-
\calL(\sseq{X}{j}{n})}.
\eeq
This observation, together with (\ref{eq:infnorm}), implies a
sufficient condition for $\eta$-mixing:
\beqn
\sum_{k=1}^\infty \phi_k
&<& \infty;
\eeqn
this certainly holds if $(\phi_k)$ admits a geometric decay, as
assumed in~\cite{samson00}.

Although $\eta$-mixing seems to be a stronger condition than
$\phi$-mixing (the latter only requires $\phi_k\to0$),
we are presently unable to obtain any nontrivial implications 
(or non-implications)
between
$\eta$-mixing 
and 
either
$\phi$-mixing
or
any of the other strong mixing conditions discussed
in~\cite{bradley}.

\subsection{Comparison between
$\nrm{\Gn}_2$
and
$\nrm{\Dn}_\infty$
}
\label{sec:DnGn}
The quantities 
$\nrm{\Gn(\P)}_2$
and
$\nrm{\Dn(\P)}_\infty$
(written here with an explicit functional dependence on the
measure $\P$)
are both 
numerical quantifiers
of the mixing properties of 
$\P$.
Because of their role in the bounds (\ref{eq:samson}) and 
(\ref{eq:dbarDn}), a smaller value for either quantity implies a
tighter deviation bound.
It turns out that 
neither 
is uniformly asymptotically tighter than the other;\hide{
We 
observe here
that neither of
$\nrm{\Gn}_2$
and
$\nrm{\Dn}_\infty$,
as an indicator
of the mixing properties of the process,
is uniformly asymptotically tighter than the other.
What we mean by
this 
precisely
is 
stated }
this statement is made precise
in 
Theorem~\ref{thm:DnGn}.
We will first need an auxiliary lemma:
\belen
\label{lem:goodbadseq}
There exists
a family of
probability spaces 
$(\X^n,\calF^n,\Pn)_{n\geq1}$ 
such that
\beqn
\label{eq:1/(n-i)}
\bar\h_{ij}(\Pn) = 1/(n-i)
\eeqn
for $1\leq i<j\leq n$.
\enlen
\begin{rem}
Since different measures are being discussed, our notation will make
explicit the functional dependence of $\bar\h_{ij}$ on the measure.
\end{rem}
\bepf
Let $\X=\set{0,1}$. 
For $1\leq k<n$,
we will call $x\in\set{0,1}^n$ a $k$-{\it good} sequence if
$x_k=x_n$ and a $k$-{\it bad} sequence otherwise. Define 
$A_n\supr k\subset \set{0,1}^n$
to be the set of the $k$-good sequences and 
$B_n\supr k=\set{0,1}^n\setminus A_n\supr k$ to be 
the bad sequences; note that $\dsabs{A_n\supr k}=\dsabs{B_n\supr k}=2^{n-1}$.
Let $\Pn\supr0$ be the uniform measure on $\set{0,1}^n$:
$$ \Pn\supr0(x) = 2^{-n},\qquad x\in\set{0,1}.$$
Now take $k=1$ and define,
for some $p_k\in[0,1/2]$,
\beqn
\label{eq:pkupd}
\Pn\supr k(x) &=& \alpha_k\Pn\supr{k-1}(x)\paren{
p_k\pred{x\in A_n\supr k} +
(1-p_k)\pred{x\in B_n\supr k}},
\eeqn
where $
\alpha_k
$ is 
the normalizing constant,
chosen so that
$\sum_{x\in\set{0,1}^n} \Pn\supr k(x)=1$.

We will say that
a probability measure $\P$ on $\set{0,1}^n$
is
$k$-{\it row homogeneous} if for all $1\leq \ell\leq k$ we have
\bit
\item[(a)]
$h_\ell(\P) \defeq \bar\h_{\ell,\ell+1}(\P) = \bar\h_{\ell,\ell+2}(\P) =\ldots = \bar\h_{\ell,n}(\P)$
\item[(b)]
$\bar\h_{ij}(\P)=0$ for $k<i<j$
\item[(c)]
$h_k$ is a continuous function of $p_k\in[0,1/2]$, with 
$h_k(0)=1$ and $h_k(1/2)=0$.
\eit
It is straightforward to verify that $\Pn\supr1$, as constructed in
(\ref{eq:pkupd}), is $1$-row homogeneous.\footnote{
The continuity of $h_k$ follows from Lemma 6.1 in \cite{kontram06}.}
Therefore,
we may choose $p_1$ 
in (\ref{eq:pkupd})
so that $h_1=1/(n-1)$. 
Iterating the formula in
(\ref{eq:pkupd})
we obtain the
sequence of measures
$\set{\Pn\supr{k}:1\leq k < n}$; each $\Pn\supr k$ is easily seen to
be $k$-row homogeneous. 
Another easily verified observation is that 
$h_{\ell}(\Pn\supr{k})=h_{\ell}(\Pn\supr{k+1})$
for all $1\leq k < n-1$ and $1\leq\ell\leq k$.
This means that we can choose the
$\set{p_k}$ 
so that $h_k(\Pn\supr{k})=1/(n-k)$ for each 
$1\leq k < n$.\hide{
we obtain the new measure 
$\Pn\supr{k+1}$; it is easy to see that
$\Pn\supr{k+1}$ is $(k+1)$-row homogeneous and that for
$1\leq\ell\leq k$, the $\ell$th row constants of $\Pn\supr{k}$ and
$\Pn\supr{k+1}$ agree. Iterating (\ref{eq:pkupd}), 
we obtain the
sequence of measures
$\Pn\supr{1},\Pn\supr{2},\ldots,\Pn\supr{n-1}$, each time choosing
$p_k$ so that $h_k=1/(n-k)$.}
The measure $\Pn\defeq \Pn\supr{n-1}$ has the desired property
(\ref{eq:1/(n-i)}).
\enpf

\bethn
\label{thm:DnGn}
There exist
families of
probability spaces 
$(\X^n,\calF^n,\Pn)_{n\geq1}$ 
such that
$R_n\to0$
and also such that
$R_n\to\infty$,
where
\beq
R_n &\defeq&
\frac{\nrm{\Gn(\Pn)}_2}{\nrm{\Dn(\Pn)}_\infty}.
\eeq
\enthn
\bepf
Recall that for an $n\times n$ real matrix $A$, its 
$\ell_\infty$ operator norm
is given by (\ref{eq:infnorm}) and its
$\ell_2$ operator norm is given by
\beq
\nrm{A}_2
=
\sup_{0\neq x\in\RR^n}\frac{\nrm{Ax}_2}{\nrm{x}_2}
=\sqrt{\lambda_{\max}(A\trn A)}
\eeq
where 
$\lambda_{\max}$ is the spectral radius.
We use the standard asymptotic ``big O'' notation, where if
$f,g:\NN\to\RR^+$, we say $f=O(g)$ if 
$\limsup_{n\to\infty}f(n)/g(n)<\infty$. The preceding relationship
between $f$ and $g$ may also be expressed as $g=\Omega(f)$. If both
$f=O(g)$ and $f=\Omega(g)$ hold, we write $f=\Theta(g)$.

For the first direction, let $\X=\set{0,1}$ and let $\Pn$ be the
measure constructed in Lemma~\ref{lem:goodbadseq},
satisfying
(\ref{eq:1/(n-i)}). For this measure, we have
$\nrm{\Dn(\Pn)}_\infty=2$ for all $n\in\NN$, so we proceed to
lower-bound $\nrm{\Gn(\Pn)}_2$. 
Letting $G_n\defeq \Gn(\P)\trn \Gn(\P)$,
an easy calculation (using 
(\ref{eq:Gndef})
and
(\ref{eq:1/(n-i)}))
gives
\beq
(G_n)_{ij} &=&
\pred{i=j} + \pred{i<j}(n-i)^{-1/2} + \pred{j<i}(n-j)^{-1/2}
+ \sum_{k=1}^{\min(i,j)-1} (n-k)\inv
\eeq
(here, $0/0\defeq0$).
Taking $x\in\RR^n$ with $x_i= i$
for $1\leq i\leq n$ and noting that
\beq
\sum_{1\leq i,j\leq n} ij\min(i,j) = \Theta(n^5),
\eeq
we 
conclude that $x\trn G_n x = \Omega(n^4)$.
Now
\beq
\nrm{x}_2 = \paren{\sum_{i=1}^n i^2}^{1/2} = \Theta(n^{3/2}),
\eeq
so
\beq
\nrm{\Gn(\P)}_2 &\geq &
\frac{\sqrt{x\trn G_n x}}{\nrm{x}_2} 
=\frac{\Omega(n^{2})}{\Theta(n^{3/2})}
\\
&=& \Omega(n^{1/2})
\eeq
and $R_n=\Omega(n^{1/2})$.

For the other direction, let $\X=\set{0,1}$ and call
$\sseq{x}{1}{n}\in\X^n$ a {\em forbidden} sequence
if $x_1\neq x_n$ and an {\em allowed} sequence otherwise.
Define the measure
$\Pn$ on $\X^n$ as vanishing on the forbidden sequences and
equiprobable on the allowed sequences:
\beqn
\label{eq:forbidden}
\Pn(\sseq{x}{1}{n}) &=& 2^{-n+1} 
\pred{x_1=x_n}.
\eeqn
For
this 
measure,
it is easy to see that
\beq
\bar\h_{ij} &=& 
\pred{i=1},
\qquad 1\leq i<j\leq n.
\eeq
This forces
$\nrm{\Dn(\Pn)}_\infty=n$
and 
\beq
(G_n)_{ij} &=& \pred{i\neq1}\pred{i=j} + 1,
\eeq
where, as before,
$G_n\defeq \Gn(\Pn)\trn \Gn(\Pn)$. To upper-bound
$\lambda_{\max}(G_n)$, 
we use a consequence of the Ger\v{s}gorin disc theorem
(\cite{hornjohnson}, 6.1.5) -- namely, that
\beq
\lambda_{\max}(G_n) &\leq&
\max_{1\leq i\leq n} \sum_{j=1}^n (G_n)_{ij}
= n+1.
\eeq
This implies $R_n = O(n^{-1/2})$.
\enpf

\begin{rem}
The last example in the proof illustrates the simple but important
point that the choice of enumeration of the random variables
$\set{X_i}$ makes a difference. Let $\pi$ be the permutation on
$\set{1,\ldots,n}$ that exchanges $2$ and $n$, leaving the other
elements fixed. 
Let $(X_i)_{1\leq i\leq n}$ be the random process on $\set{0,1}^n$ defined in
(\ref{eq:forbidden}) and define
process 
$Y=\pi(X)$ by
$Y_i = X_{\pi(i)}$, $1\leq i\leq n$.
It is easily verified that
$\nrm{\Dn(Y)}_\infty=2$ 
while we saw above that 
$\nrm{\Dn(X)}_\infty=n$. Thus if $f:\set{0,1}^n\to\RR$ is invariant
under permutations and $\xi_1,\xi_2\in\RR$ are random variables
defined by
$\xi_1=f(X)$,
$\xi_2=f(\pi(X))$,
we have $\xi_1=\xi_2$ with probability $1$, yet our technique proves
much tighter concentration for $\xi_2$ than for $\xi_1$.
Of course, knowing this special relationship between $\xi_1$ and
$\xi_2$, we can deduce a corresponding concentration result for
$\xi_1$; what is crucial is that the concentration for $\xi_1$ is
obtained by re-indexing the random variables.
\end{rem}

\begin{rem}
Note that for the first direction in the proof of
Theorem~\ref{thm:DnGn}, we constructed a sequence of measures $\Pn$
such that
$\nrm{\Dn(\Pn)}_\infty=2$ is bounded while
$\nrm{\Gn(\Pn)}_2=\Omega(n^{1/2})$. Is there a sequence of measures
for which
$\nrm{\Gn(\Pn)}_2$ is bounded and
$\nrm{\Dn(\Pn)}_\infty$ unbounded? We conjecture that such a sequence
of measures indeed exists, but leave its construction for future investigation.
\hide{
On a coarser scale, it seems reasonable to conjecture that 
each
of
$\nrm{\Gn(\Pn)}_2$ and $\nrm{\Dn(\Pn)}_\infty$
will be bounded
if and only if the other one is;
we leave this
question open for future investigation.
}
\end{rem}
\begin{rem}
In Lemma~\ref{lem:goodbadseq}, we constructed a sequence of measures
$\Pn$ so that $\Dn(\Pn)$ has a specific form. An obvious constraint on
the form of
$\Dn$ is 
\\
(*) $0\leq \bar\h_{ij}\leq 1$,\\
and the constraint\\
(**) $\bar\h_{i,j}\geq\bar\h_{i,j+1}$, for $1\leq i<j<n$\\
is easily seen to hold for all measures $\Pn$ on $\X^n$. Do (*) and
(**) completely specify the set of the possible $\Dn(\Pn)$ -- or
are there other constraints that all such matrices must satisfy? We
are inclined to conjecture the former, but leave this question open
for now.
\end{rem}

\hide{
On the other hand, if we disregard the asymptotic rates and restrict
ourselves to the cruder question of boundedness, we get a positive
result:
\bethn
For any
family of
probability spaces 
$(\X^n,\calF^n,\Pn)_{n\in\NN}$,
\beq
\sup_{n\in\NN}\nrm{\Gn(\Pn)}_2 <\infty
&\Longleftrightarrow&
\sup_{n\in\NN}\nrm{\Dn(\Pn)}_\infty <\infty
\eeq
\enthn
\bepf
Let $G_n$ be the upper-triangular matrix defined by
$G_n \defeq \Gn(\Pn)$. Now the 
$n\times n$ matrices form a finite ($n^2$) dimensional vector space,
so the
$\ell_2$ and $\ell_\infty$ operator norms are equivalent. That means
that there exist constants $c_n,c'_n>0$, depending only on $n$, such
that
\beq
c_n\nrm{G_n}_2 \;\leq\;
\nrm{G_n}_\infty \;\leq\;
c'_n\nrm{G_n}_2.
\eeq
Let $D_n = \Dn(\Pn)$; thus $(H_n)_{ij}^2 = (D_n)_{ij}$.

\enpf
}

\section{$\Psi$-dominance}
\label{sec:psidom}
Having dealt with the ``analytic'' mixing condition on $X$ in \S\ref{sec:hmix}, we
now turn to the 
geometry of $(\X^n,\rho_n)$.

We say that the family of metric measure spaces 
$(\X^n,\rho_n,\mu^n)_{n\geq1}$ is
{\em consistent} if
\bit
\item[(\one)] 
the metrics $\{\rho_n\}$
satisfy, for all 
$1\leq i\leq n$
and
$\sseq{x}{1}{n},\sseq{y}{1}{n}\in\X^n$,
\beq
\rho_n(\sseq{x}{1}{n},\sseq{y}{1}{n})=
\rho_{n-1}(
\scat{\sseq{x}{1}{i-1}\sd\sseq{x}{i+1}{n}},
\scat{\sseq{y}{1}{i-1}\sd\sseq{y}{i+1}{n}}),
\eeq
whenever $x_i=y_i$
\item[(\two)]
for each $n\geq1$,
$\mu^n$ is a positive product measure on the Borel
$\sigma$-algebra induced by $\rho_n$.
\hide{
\item[(\thr)]
for $n\geq1$,
$\diam_{\rho_n}(\X^n)\leq n$.
}
\eit
\begin{rem}
Condition (\one) implies that the topology $\tau^n$ induced by $\rho_n$ on
$\X^n$ is the product topology 
$\tau^n = \tau\otimes\tau\otimes\ldots\otimes\tau$,
where $\tau$ is the topology induced on $\X$ by $\rho_1$.
Likewise, $\mu$ 
is a positive measure on the Borel $\sigma$-algebra generated by
$(\X,\rho_1)$
and
$\mu^n=\mu\otimes\mu\otimes\ldots\otimes\mu$
is the corresponding product measure on the product 
$\sigma$-algebra. 
\end{rem}

A quantitative notion of continuity
is the Lipschitz condition, which is defined with respect to some
metric $\rho_n$ on $\X^n$. 
Define $\Lips(\X^n,\rho_n)$ to be the set of all
$f:\X^n\to[0,\diam_{\rho_n}(\X^n)]$ such that
\beqn
\label{eq:Lips}
\sup_{x\neq y\in\X^n}
\frac{
\abs{f(x)-f(y)}
}{\rho_n(x,y)}
 &\leq& 
1
\eeqn
(any such function is continuous and therefore measurable).
\begin{rem}
\label{rem:[0,n]}
Since the Lipschitz condition implies $\diam f(\X^n)\leq \diam\X^n$
and the functionals 
$V_i$ 
and $\hat V_i$ 
(defined in
(\ref{eq:videf})
and
(\ref{eq:vhat}),
respectively)
are translation-invariant (in the sense that $V_i(f;y)=V_i(f+a;y)$ for
all $a\in\RR$), there is no loss of generality in restricting the
range of $f$ to $[0,\diam\X^n]$.
\end{rem}

\hide{
As in \S\ref{sec:simplemartbd}, we assume a 
positive
Borel carrying measure
$\mu$ on $\X^n$ such that $\P\abscont\mu$ with density 
$d\P(x)=p(x)d\mu(x)$.
}

Let $F_n=L_1(\X^n,\mu^n)$ 
\hide{
be the space of all $\mu$-integrable functions
$f:\X^n\to\RR$. 
Let $\Phi_n\subset F_n$ be the functions $g\in F_n$
with range in $[0,n]$
satisfying the Lipschitz condition
\beqn
\label{eq:lipcond}
 \abs{g(x)-g(y)} &\leq& \nrm{x-y}_1
\eeqn
for all $x,y\in[0,1]^n$,
with $\nrm{\cdot}_1$ being the $\ell_1$ norm on $\RR^n$. 
}
and
equip $F_n$
with the inner product 
\beqn
\label{eq:ipdef}
 \iprod{f}{g} &=& \int_{\X^n} f(x)g(x) d\mu^n(x).
\eeqn
Since $f,g\in F_n$ might not be in $L_2(\X^n,\mu^n)$, the expression in
(\ref{eq:ipdef})
in general
might not be finite. However, for $g\in\Lips(\X^n,\rho_n)$, we
have
\beqn
\label{eq:phitriv}
\abs{\iprod{f}{g}} &\leq& \diam_{\rho_n}(\X^n) \nrm{f}_{L_1(\mu^n)}
\eeqn
(the motivation for 
bounding
$\iprod{f}{g}$
comes from
(\ref{eq:linfunc})).

Define the {\em marginal projection} operator $\pi:F_{n}\to F_{n-1}$ as
follows. If 
$f:\X^n\to\RR$ then
$(\pi f):\X^{n-1}\to\RR$ is given by
\beqn
\label{eq:pidef}
(\pi f)
(x_2,\ldots,x_n) &\defeq& \int_\X f(x_1,x_2,\ldots,x_n) d\mu(x_1).
\eeqn
Note that by Fubini's theorem (Thm. 8.8(c) in \cite{rudin}), 
$\pi f \in L_1(\X^{n-1},\mu^{n-1})$.
Define the functional $\psin{n}:F_n\to\RR$
recursively: $\psin{0}\defeq 0$ and
\beqn
\label{eq:contpsidef}
 \psin{n}(f) &\defeq& \int_{\X^n} \pl{f(x)}d\mu^n(x) + \psin{n-1}(\pi f) 
\eeqn
for $n\geq1$. 
The latter is finite since
\beqn
\label{eq:psitriv}
\psin{n}(f) &\leq& n\nrm{f}_{L_1(\mu)},
\eeqn
as shown in Theorem~\ref{thm:psinorm} below.

We say that the family of metric spaces $(\X^n,\rho_n)_{n\geq1}$ is
{\em $\Psi$-dominated} 
with respect to 
a positive Borel measure $\mu$ on $\X$
if 
$(\X^n,\rho_n,\mu^n)_{n\geq1}$ is consistent in the sense of
(\one) and (\two) above, and
the inequality
\beqn
\label{eq:psidom}
\sup_{g\in\Lips(\X^n,\rho_n)} \iprod{f}{g} &\leq& \psin{n}(f)
\eeqn
holds
for all $f\in L_1(\X^n,\mu^n)$.

\bethn
\label{thm:taurhodom}
Suppose 
$(\X^n,\rho_n)_{n\geq1}$ is a
$\Psi$-dominated
family of metric spaces
with respect to some 
(positive Borel)
measure $\mu$
and
$(\X^n,\tau_n)_{n\geq1}$ is another family of metric spaces,
with $\tau_n$ dominated by $\rho_n$, 
in the sense
that
\beqn
\label{eq:taurho}
\tau_n(x,y) &\leq& \rho_n(x,y),
\qquad x,y\in\X^n
\eeqn
for all $n\geq1$.
Then $(\X^n,\tau_n)_{n\geq1}$ is also
$\Psi$-dominated 
with respect to 
$\mu$.
\enthn
\bepf
By (\ref{eq:taurho}), we have
\beq
\Lips(\X^n,\tau_n)
&
\subset
&
\Lips(\X^n,\rho_n),
\eeq
which in turn implies
\beq
\sup_{g\in \Lips(\X^n,\tau_n)}
\abs{\iprod{f}{g}}
&\leq&
\sup_{g\in \Lips(\X^n,\rho_n)}
\abs{\iprod{f}{g}}
\leq \psin{n}(f).
\eeq

\enpf

We are about to define two functionals on 
$F_n=L_1(\X^n,\mu^n)$.
Although we use the norm
notation, none of the results we prove actually rely on the norm
properties of 
$\phinorm{\cdot}$ and 
$\psinorm{\cdot}$, 
and
so we defer a discussion of these
do the Appendix. The punchline is that under appropriate conditions both
are valid norms;
$\psinorm{\cdot}$ is (topologically) equivalent to $\nrm{\cdot}_{L_1}$ while
$\phinorm{\cdot}$ is in general weaker.

The two norms are defined as
\beqn
\label{eq:phinorm}
\phinorm{f} &=& \sup_{g\in \Lips(\X^n,\rho_n)} \abs{\iprod{f}{g}}
\eeqn
and
\beqn
\label{eq:psinorm}
 \psinorm{f} &=& \max_{s=\pm1}\psin{n}(s f);
\eeqn
note that 
(\ref{eq:psidom}) 
is equivalent to
the condition that
$\phinorm{f}\leq\psinorm{f}$ for all $f\in F_n$.
We refer to the norms in
(\ref{eq:phinorm})  and (\ref{eq:psinorm})
as $\Phi$-norm and $\Psi$-norm, respectively;
notice that both depend on the measure $\mu$ and $\Phi$-norm also
depends on the metric.

\section{Main result:
$\eta$-mixing
with
$\Psi$-dominance
imply normal L\'evy
family}
\label{sec:main}
\bethn
\label{thm:main}
Let $(\X^k,\rho_k,\P)_{1\leq k\leq n}$ be a 
$\Psi$-dominated
family of
metric probability spaces
with respect to 
a
positive Borel measure $\mu$,
where
$\P\abscont\mu^n$.
Then, for any 
Lipschitz
(with respect to $\rho_n$)
$f:\X^n\to\RR$
we have
\beq
\pr{\abs{f-\E f}>
t
} &\leq&
2\exp\paren{-\frac{t^2}{2
n
\Lip{f}^2
\nrm{\Dn}_\infty^2}}
\eeq
for all $t>0$,
where 
$\Dn$
is defined in (\ref{eq:Ddef}).
\hide{
For $f\in\Lips(\X^n,\rho_n)$
let
$V_i(f;\cdot)$ 
be 
the martingale difference
defined in
(\ref{eq:videf}), for $1\leq i\leq n$.
Then we have
\beq
\nrm{V_i(f;\cdot)}_{L_\infty(\P)}
&\leq&
1+\sum_{j=i+1}^n \bar\h_{ij}.
\eeq
}
\enthn
\begin{rem}
A version of this result is proved in 
Theorem 5.1
of~\cite{kontram06},
for the special case of the counting measure 
on a finite set
$\X^n$, where $\rho$ is the Hamming metric. 
Note that if we require $\Lip{f}\leq1$ with respect to the normalized
metric
$\bar\rho_n=\oo n \rho_n$, we get
\beqn
\label{eq:dbarDn}
\pr{\abs{f-\E f}>
t
} &\leq&
2\exp\paren{-\frac{nt^2}{2
\nrm{\Dn}_\infty^2}};
\eeqn
for $\eta$-mixing measures $\P$
 (see (\ref{eq:HP})),
this implies
$\pr{\abs{f-\E f}>t} \leq 2\exp(-nt^2/2 \bar H_{\P})$, meaning that
the
$(\X^n,\bar\rho_n,\P)$ form a normal L\'evy family.

We will use the same conventions regarding the density
$d\P(x) = p(x)d\mu^n(x)$ as in
\S\ref{sec:simplemartbd}.
\end{rem}
\bepf
The claim will follow 
via (\ref{eq:genconc}),
by 
proving
the
bound
\beqn
\label{eq:VD}
\tsnrm{V_i(f;\cdot)}_{L_\infty(\P)} &\leq& \Lip{f}\tsnrm{\Dn}_\infty
\eeqn
on
the martingale difference $V_i(f;\cdot)$.
Since
$\tsnrm{V_i(f;\cdot)}_{L_\infty}$
and
$\Lip{f}$
are both homogeneous functionals of $f$ (in the sense of
$T(af)=|a|T(f)$ for $a\in\RR$), there is no loss of generality in
taking $\Lip{f}=1$.

Lemma~\ref{lem:vhat} shows that
it suffices to bound $\tsnrm{\hat V_i(f;\cdot)}_{L_\infty}$, and
from (\ref{eq:linfunc}), we have
\beqn
\label{eq:mgaleiprod}
\hat V_i(f;\yi,w_i,w_i') \;=\;
\int_{\X^n} f(x) \hat g(x)d\mu^n(x)
\;=\;
\iprod{f}{\hat g},
\eeqn
where
\beqn
\label{eq:ghat}
\hat g(x) &=&
\pred{\sii= \scat{\yi\sd w_i}}
p(\xin\gn \scat{\yi\sd w_i})
-
\pred{\sii =\scat{\yi\sd w_i'}}
p(\xin\gn \scat{\yi\sd w_i'}).
\eeqn
Let $1\leq i<j\leq n$ and $y\in\X^{i-1},w,w'\in\X$ be fixed.
For $1\leq k\leq n$,
let
$F_k=L_1(\X^k,\mu^k)$ and recall the definition 
(\ref{eq:pidef})
of the projection
operator $\pi:F_{k}\to F_{k-1}$.
Put $N=n-i+1$ and 
for $y\in\X^{i-1}$
define 
the operator 
$T_y:F_n\to F_N$
by
\beq
(T_y f)(x) &\defeq& f(yx)
\eeq
for each $x\in\X^N$.
Observe 
that
(\ref{eq:ghat}) implies
\hide{
$f_{iy},\hat g_{iy}\in F_N$ by
\beq
     f_{iy}(z) =      f([y\,z]), \qquad
\hat g_{iy}(z) = \hat g([y\,z])
\eeq
for all $z\in\X^N$,
noting that}
\beqn
\label{eq:fgiy}
\iprod{f}{\hat g} &=&
\iprod{T_yf}{T_y\hat g}.
\eeqn
By Remark~\ref{rem:[0,n]}, we may take $f\in\Lips(\X^n,\rho_n)$, and
therefore 
(by 
the consistency of the metrics, in the sense of \S\ref{sec:psidom}),
$T_yf\in\Lips(\X^N,\rho_N)$.

Let $\hat g\supr{N}\defeq T_y\hat g$ and for
$\ell=N,N-1,\ldots,2$, define
\beq
\hat g\supr{\ell-1} &=& \pi \hat g\supr{\ell};
\eeq
note that $\hat g\supr{\ell}\in F_\ell$.

A direct calculation (using
the 
Radon-Nikod\'ym theorem) gives
\beq
\hat g\supr{n-j+1}(x) &=&
p(\sseq{X}{j}{n}=x\gn \sseq{X}{1}{i}=\scat{y\sd w}) -
p(\sseq{X}{j}{n}=x\gn \sseq{X}{1}{i}=\scat{y\sd w'})
\eeq
for all $x\in\X^{n-j+1}$. It follows via
(\ref{eq:taubal}) that
\beq
\h_{ij}(y,w,w')
&=&
\int_{\X^{n-j+1}} \pl{\hat g\supr{n-j+1}(x)} d\mu^{n-j+1}(x).
\eeq
Since the measure $d\nu = \hat g\supr{n-i+1}(x) d\mu^{n-i+1}(x)$ is the difference
of two probability measures, we have $\TV{\nu}\leq1$. 
Thus the
definition of
the $\psin{n-i+1}$ functional 
(acting on $F_{n-i+1}$)
yields
\beq
\psinorm{T_y\hat g} 
&\leq& 1+\sum_{j=i+1}^n \h_{ij}(y,w,w') \\
&\leq& 1+\sum_{j=i+1}^n \bar\h_{ij} \qquad\mbox{$\P$-almost surely}\\
&\leq& \nrm{\Dn}_\infty.
\eeq
Putting together
(\ref{eq:psidom}),
(\ref{eq:mgaleiprod})
and
(\ref{eq:fgiy}),
we obtain the desired bound in
(\ref{eq:VD}).
\enpf

\section{
Applications
}
\label{sec:pdex}
\subsection{$(\NN^n,d_\ham)$ is $\Psi$-dominated}
A core result in~\cite{kontram06} (Theorem 4.8)
effectively
established the $\Psi$-dominance of
$(\X^n,d_\ham)$ for finite $\X$.
For the countable case, verifying
consistency (in the sense of \S\ref{sec:psidom})
is trivial. Let
$\X=\NN$, $\mu$ be the counting measure on $\X^n$,
$f\in 
\ell_1(\X^n)\equiv L_1(\X^n,\mu)
$ and $g\in\Lips(\X^n,d_\ham)$.
For $m\geq1$,
let
$\X_m = \{k\in\X: k\leq m\}$
and define
the $m$-{\em truncation} of $f$ to be
the following function in
$\ell_1(\X^n)$:
\beq
f_{m}(x) 
&=&
\pred{x\in\X_m^n}f(x).
\eeq
Then we have, by~\cite{kontram06}, Theorem 4.8,
\beq
\iprod{f_m}{g} &\leq& \psin{n}(f_m)
\eeq
for all $m\geq1$, and $\lim_{m\to\infty}f_m(x)=f(x)$ for all
$x\in\X^n$. 
Let $h_m(x) = f_m(x)g(x)$ and note that
$|h_m(x)| \leq n|f(x)|$, the latter in 
$\ell_1(\X^n)$. 
Thus by
Lebesgue's Dominated Convergence theorem, we have 
$\iprod{f_m}{g}\to\iprod{f}{g}$.
A similar dominated convergence argument
shows that
$\psin{n}(f_m)\to\psin{n}(f)$, which proves the
$\Psi$-dominance of
$(\NN^n,d_\ham)$.


\subsection{$([0,1]^n,\nrm{\cdot}_1)$ is $\Psi$-dominated}
\newcommand{\Phin}{\Lips([0,1]^n,\nrm{\cdot}_1)}
\newcommand{\Phimn}{\Lips(\X_m^n,d_m)}
\newcommand{\Phiham}{\Lips(\X_m^n,d_\ham)}
Since verifying 
consistency 
is trivial,
it remains to prove
\bethn
\label{thm:l1psidom}
Let $\mu$ be the Lebesgue measure on $[0,1]$ and
$\rho_{n}(x,y)=\nrm{x-y}_1$, for $x,y\in[0,1]^n$. 
Then we have
\beqn
\label{eq:phipsiL1}
 \phinorm{f} &\leq& \psinorm{f}
\eeqn
for all $f\in L_1([0,1]^n,\mu^n)$.
\enthn
\bepf
\hide{
First observe that from (\ref{eq:phitriv}) and
(\ref{eq:psitriv}), we have
\beqn
\label{eq:phipsiL1} 
(\phinorm{f}<\infty) \iff
(\psinorm{f}<\infty) \iff
(\nrm{f}_1<\infty)
\eeqn
for all measurable $f:[0,1]^n\to\RR$.
Lemma~\ref{lem:phipsitriv} implies
that the norms 
$\phinorm{\cdot}$, 
$\psinorm{\cdot}$,
and
$\nrm{\cdot}_1$
are equivalent in the sense of defining the same class of Cauchy
sequences (and hence topology) on $F_n$.
}
Let $F_n=L_1([0,1]^n,\mu^n)$
and
$C_n\subset F_n$ be the class of continuous 
functions.
It follows from Theorem 3.14
of~\cite{rudin} that $C_n$ is dense in $F_n$,
in the topology induced by $\nrm{\cdot}_{L_1}$.
This
implies
that
for any $f\in
F_n$ and $\eps>0$, there is a $g\in C_n$ such that
$\nrm{f-g}_{L_1}<\eps/n$
and therefore
(via 
(\ref{eq:phitriv}) and
(\ref{eq:psitriv})),
$$ 
\phinorm{f-g}<\eps
\qquad\text{and}\qquad
\psinorm{f-g}
< \eps
,$$
so it suffices to prove (\ref{eq:phipsiL1}) for $f\in C_n$.

For 
$m>1$,
let $\X_m=\{k\in\NN:0\leq k< m\}$. 
Define the {\em grid map}
$\gamma_m:[0,1]^n\to\X_m^n$ by
\beq
[\gamma_m(x)]_i &=& 
\max\{k\in\X_m : k/m\leq x_i\}
\eeq
for $x\in[0,1]^n$ and $1\leq i \leq n$;
$x
$ is called an $m$-{\em grid point} if each coordinate $x_i$ is
of the form $x_i=s/m$, for some $s\in\X_m$.
 We say that $g\in F_n$ is a
{\em grid-constant function} if there is an $m>1$ such that $g(x)=g(y)$
whenever $\gamma_m(x)=\gamma_m(y)$; thus a grid-constant 
function is constant on the grid cells.
Let $G_n\subset
F_n$ be the class of 
grid-constant
functions. It is easy to see that $G_n$ is
dense in $C_n$. 
Indeed,
for $f\in C_n$ and $\eps>0$, there is a $\delta>0$
such that $\omega_f(\delta)<\eps$, where $\omega_f$ is the 
$\ell_\infty$
modulus of
continuity of $f$. Taking $m=\ceil{1/\delta}$ and $g\in G_n$ to be
such that it agrees with $f$ on the $m$-grid points, we have
$\nrm{f-g}_{L_1([0,1]^n)}\leq\nrm{f-g}_{L_\infty([0,1]^n)}<\eps$. 
Thus we need only prove (\ref{eq:phipsiL1}) for $f\in G_n$.

Define the metric $d_m$ on $\X_m$:
$$ d_m(z,z') = \frac{\abs{z-z'}}{m-1}$$
and extend it to $\X_m^n$:
$$ d_m(z,z') = \sum_{i=1}^n d_m(z_i,z_i').$$
Let $D_n\subset G_n$ consist of those functions 
$g:[0,1]^n\to[0,n]$ 
for which there is
an $m>1$ such that
$$ \abs{g(x)-g(y)} \leq d_m(\gamma_m(x),\gamma_m(y))$$
for all $x,y\in[0,1]^n$. The argument used above shows that $D_n$ is
dense in $\Phin$, and so it suffices to bound
$\sup_{g\in D_n} \iprod{f}{g}$
for $f\in G_n$. 

Fix $f\in G_n$, $g\in D_n$, and let $m>1$
be such that $f$ and $g$ are $m$-grid-constant functions.
Let 
$\bar\k,\bar\f:\X_m^n\to\RR$ 
be such that 
$\bar\k(\gamma_m(x))=f(x)$
and
$\bar\f(\gamma_m(x))=g(x)$
for all $x\in[0,1]^n$.
Then
\beq
\iprod{f}{g} &=& \paren{\oo m}^n \sum_{z\in\X_m^n} \bar\k(z)\bar\f(z)
\eeq
and
\beq
\psin{n}(f) &=& \paren{\oo m}^n \hat\Psi_n(\bar\k),
\eeq
where $\hat\Psi_n$ 
is 
$\Psi_n$ 
computed using the counting measure on $\X_m$.

Define $\Phimn$ and $\Phiham$ in accordance with
(\ref{eq:Lips})
and note that
$\bar\f\in\Phimn$.\hide{
Define $\hat\Phi_{n,m}$ to be the set of all $\f:\X_m^n\to[0,n]$ such that
\beq
\abs{\f(z)-\f(z')}&\leq& d_m(z,z'),\qquad z,z'\in\X_m^n
\eeq
(note that $\bar\f\in\hat\Phi_{n,m}$)
and let $\hat\Phi_{n,\ham}$ consist of those $\f:\X_m^n\to[0,n]$ for which
\beq
\abs{\f(z)-\f(z')}&\leq& d_{\ham}(z,z'),\qquad z,z'\in\X_m^n
\eeq
where $d_{\ham}$ is the Hamming metric on $\X_m^n$.}
We claim that 
$\Phimn\subset\Phiham$; 
this holds because
$d_m(z,z')\leq d_{\ham}(z,z')$. 
Theorem 4.8
in~\cite{kontram06}
states that
for all $\k:\X_m^n\to\RR$,
\beq
\sup_{\f\in \Phiham} \sum_{z\in\X_m^n}\k(z)\f(z)
&\leq&
\hat\Psi_n(\k).
\eeq
This
implies $\iprod{f}{g}\leq\psin{n}(f)$ and completes the proof.
\enpf

\begin{rem}
One might be tempted to take a shortcut to this result by showing 
directly
that
$([0,1]^n,d_\ham)$ is $\Psi$-dominated
and then applying Theorem~\ref{thm:taurhodom} to $d_\ham$ and $\nrm{\cdot}_1$.
The problem with this approach
is that $d_\ham$ induces the discrete topology on $[0,1]^n$, whose
open sets are not necessarily Lebesgue measurable.
\end{rem}

\subsection{$([0,1]^n,
\nrm{\cdot}_p
)$ is $\Psi$-dominated}
\label{sec:pdom}
Recall that
for any $1< p\leq\infty$ and any
$x\in\RR^n$, we have
\beqn
\label{eq:l1p}
\nrm{x}_p \;\leq\;
\nrm{x}_1
\;\leq\;
n^{1/p'}\nrm{x}_p,
\eeqn
where $1/p+1/p'=1$. 
The 
first 
inequality 
holds
because the convex
function $x\mapsto\nrm{x}_p$ is maximized on the extreme points (corners)
of the convex polytope $\set{x\in\RR^n:\nrm{x}_1=1}$.
The second inequality is checked by applying H\"older's inequality to
$\sum x_i y_i$, with $y\equiv1$. Both 
are tight.
Furthermore, all the $\ell_p$ norms induce the same topology on
$\RR^n$, whose Borel sets are Lebesgue measurable.
Thus, in light of 
Theorem~\ref{thm:taurhodom},\hide{
(\ref{eq:l1p}), we have
\beq
f\in\Lips([0,1]^n,\nrm{\cdot}_p)
&\implies&
f\in\Lips([0,1]^n,\nrm{\cdot}_1),
\eeq
which in turn implies
\beq
\sup_{g\in \Lips([0,1]^n,\nrm{\cdot}_p)} \abs{\iprod{f}{g}}
&\leq&
\sup_{g\in \Lips([0,1]^n,\nrm{\cdot}_1)} \abs{\iprod{f}{g}}.
\eeq
Thus}
the $\Psi$-dominance (with respect to the Lebesgue measure,
see Theorem~\ref{thm:l1psidom}) of
$\nrm{\cdot}_1$
implies the $\Psi$-dominance of $\nrm{\cdot}_p$.

\hide{
Applying Theorem~\ref{thm:main} to $\rho(x,y)=\nrm{x-y}_1$ on $[0,1]^n$,
we have for any measurable $f:[0,1]^n\to\RR$
\beq
\pr{\abs{f-\E f}>t} &\leq&
2\exp\paren{-\frac{t^2}{2 n \Lipp{f}{1}^2
\nrm{\Dn}_\infty^2}}
\eeq
where $\Lipp{f}{p}$ is 
the Lipschitz constant
with respect to $\nrm{\cdot}_p$, $1\leq p\leq\infty$.
}

\subsection{Converting between Samson's bound and
  Theorem~\ref{thm:main}}
\label{sec:convert}
Let us attempt a rough comparison between the results obtained here
and the main result of Samson's 2000 paper~\cite{samson00}. In light
of Theorem~\ref{thm:DnGn}, a uniform comparison between our mixing
coefficient $\nrm{\Dn}_\infty$ and Samson's $\nrm{\Gn}_2$ is not
possible. However, assume for simplicity that for a given random
process $X$ on $[0,1]^n$, the two quantities are of the same order of
magnitude. 
For example, for the case of contracting Markov chains with Doeblin
coefficient $\tha<1$, we have
\beq
\nrm{\Dn}_\infty\leq\oo{1-\tha},
\qquad
\nrm{\Gn}_2\leq\oo{1-\tha^{1/2}}
\eeq
(as computed in~\cite{kontram06} and~\cite{samson00}, respectively).

Throughout this discussion, we will take $\X=[0,1]$
and $\mu$ to be the Lebesgue measure. For 
$f:\RR^n\to\RR$, we define
$\Lipp{f}{p}$ to be
the (smallest) Lipschitz constant of $f$
with respect to 
the metric
$d(x,y)=\nrm{x-y}_p$, 
where
$1\leq p\leq\infty$.

Suppose $f:[0,1]^n\to\RR$ has $\Lipp{f}{2}\leq1$.
Samson gives the deviation inequality
\beq
\pr{\abs{f-\mexp f}>t} &\leq& 2\exp\paren{-\frac{t^2}{2\nrm{\Gn}_2^2}}
\eeq
with the additional requirement that $f$ be convex.
By (\ref{eq:l1p}) we have $\Lipp{f}{1}\leq1$
and by Theorem~\ref{thm:l1psidom}, the $\ell_1$ metric is
$\Psi$-dominated. Thus, Theorem~\ref{thm:main} applies:
\hide{
Alternatively, we may apply
Theorems
\ref{thm:main}
and
\ref{thm:phipsiL1},
as well as
the observation in \S\ref{sec:pdom}, 
to
get}
\beqn
\label{eq:thmainapp}
\pr{\abs{f-\mexp f}>t\sqrt{n}} &\leq& 2\exp\paren{-\frac{t^2}{2\nrm{\Dn}_\infty^2}}
\eeqn
for any $f:[0,1]^n\to\RR$
with
$\Lipp{f}{2}\leq1$ (convexity is not required).

To convert from the bound in Theorem~\ref{thm:main} to Samson's bound,
we start with a convex
$f:[0,1]^n\to\RR$, 
having $\Lipp{f}{1}\leq1$. By (\ref{eq:l1p}), this means that
$\Lipp{f}{2}\leq\sqrt n$, or equivalently,
$\Lipp{n^{-1/2}f}{2}\leq1$.
Applying Samson's bound to $n^{-1/2}f$, we get
\beqn
\label{eq:samsonapp}
\pr{\abs{f-\mexp f}>t\sqrt{n}} &\leq& 2\exp\paren{-\frac{t^2}{2\nrm{\Gn}_2^2}},
\eeqn
while the bound provided by Theorem~\ref{thm:main} remains as stated
in (\ref{eq:thmainapp}).

We stress that the factor of $\sqrt n$ in 
(\ref{eq:thmainapp}) and 
(\ref{eq:samsonapp}) appears in the two bounds for rather different
reasons. In 
(\ref{eq:thmainapp}),
it is simply another way of stating
Theorems~\ref{thm:main}
and
\ref{thm:l1psidom}
for $\Lipp{f}{1}\leq 1$; namely,
$\pr{\abs{f-\mexp f}>t} \leq 2\exp(-t^2/2n \nrm{\Dn}_\infty^2)$.
In (\ref{eq:samsonapp}), the $\sqrt n$ was the ``conversion cost''
between the $\ell_1$ and the $\ell_2$ metrics.


\section{Discussion}
\label{sec:disc}
We have provided a general framework for proving measure concentration
results in various metric spaces. A useful feature of our treatment is
its modularity: since the geometric properties of the metric
($\Psi$-dominance) have been decoupled from the analytic properties of
the measure ($\eta$-mixing), Theorem~\ref{thm:main} actually gives
rise to a family of measure concentration results.

While the bounds stated in terms of $\Dn$ are not 
directly
comparable to the ones in terms of $\Gn$, we provide some discussion and
intuition in \S\ref{sec:DnGn} and \S\ref{sec:convert}. The rough
summary is that neither gives asymptotically tighter bounds 
than the other
uniformly over
all
processes, and that the former is most suitable for the $\ell_1$ metric
while the latter works best with $\ell_2$ (though both are applicable
to general $\ell_p$ metrics; see 
\S\ref{sec:pdom} and
\S\ref{sec:convert}). Samson's deviation inequality requires that $f$
be convex while ours does not; we also note that the
$\ell_\infty$ operator norm 
$\nrm{\Dn}_\infty$
is often simpler to estimate than the
spectral norm
$\nrm{\Gn}_2$.

Comparisons aside, we have offered a new approach for studying the
concentration of measure phenomenon and are hopeful that it will find
interesting applications in future work.

\begin{appendix}

\section{
Norm properties of 
$\phinorm{\cdot}$
and
$\psinorm{\cdot}$
}
\label{sec:norms}
It was proved in~\cite{kontram06} that
$\phinorm{\cdot}$
and
$\psinorm{\cdot}$
are valid norms when $\X$ is finite. We now do this in a significantly
more
general setting, and examine the strength of the toplogies induced by
these norms.

\bethn
\label{thm:psinorm}
Let $F_n=L_1(\X^n,\mu^n)$ for some positive Borel measure $\mu$.
Then
\bit
\item[(a)] $\psinorm{\cdot}$ is a vector-space norm on $F_n$
\item[(b)] for all $f\in F_n$,
\beq
          {\ts\oo2}\nrm{f}_{L_1} \;\leq \; \psinorm{f}
\;\leq \; n        \nrm{f}_{L_1} .
\eeq
\eit
\enthn
\bepf
We prove (b) first. Since
\beq
\tsnrm{f}_{L_1} &=& \tsnrm{\pl{f}}_{L_1} + \tsnrm{\pl{-f}}_{L_1},
\eeq
we have that 
$\psinorm{f}$
(defined in 
(\ref{eq:contpsidef})
and
(\ref{eq:psinorm})) 
is the sum of
$n$ terms, each one at most $\nrm{f}_{L_1}$ and the first one at least
${\ts\oo2}\nrm{f}_{L_1}$; this proves (b).

To prove (a) we check the norm axioms:

{\em Positivity}: 
It is obvious that $\psinorm{f}\geq0$ and
(b) shows that $\psinorm{f}=0$ and iff $f=0$ a.e. $[\mu]$.

{\em Homogeneity}: It is immediate from 
(\ref{eq:contpsidef})
that
$\psin{n}(af) = a\psin{n}(f)$
for 
$a\geq0$.
From (\ref{eq:psinorm}) we have 
$\psinorm{f} = \psinorm{-f}$. 
Together these imply
$\psinorm{af} = |a|\psinorm{f}$. 

{\em Subadditivity}: 
It follows from the subadditivity of the function $h(z)=\pl{z}$ and
additivity of integration that
$\psinorm{f+g} \leq \psinorm{f} + \psinorm{g}$.
\enpf

\bethn
\label{thm:phisemi}
Let $F_n=L_1(\X^n,\mu)$ for some 
metric measure
space $(\X^n,\rho,\mu^n)$. Then $\phinorm{\cdot}$ is a 
seminorm
on
$F_n$. 
\enthn
\bepf

{\em Nonnegativity}:
$\phinorm{f}\geq0$
is
obvious
from the definition (\ref{eq:phinorm}).

{\em Homogeneity}: It is clear 
from the definition
that 
$\phinorm{af}=|a|\phinorm{f}$
for any
$a\in\RR$.

{\em Subadditivity}:
$\phinorm{f+g} \leq \phinorm{f} + \phinorm{g}$
follows from the linearity of
$\iprod{\cdot}{\cdot}$ and the triangle inequality for $\abs{\cdot}$.
\enpf

Under mild conditions
on the 
metric 
measure
space
$(\X^n,\rho,\mu^n)$,
$\phinorm{\cdot}$ is a
genuine norm. 
We will use
the
topological notion of
local compactness
(meaning that every point has a neighborhood with compact closure).
We
also require some regularity conditions on the measure $\mu^n$:
\bit
\item[(a)] $\mu^n(K)<\infty$ for every compact set $K\subset\X^n$
\item[(b)] for every Borel $E\subset\X^n$, we have
\beq
\mu^n(E) = \inf\set{\mu^n(V): E\subset V,~V~\mbox{open}}
\eeq
\item[(c)] if $E\subset\X^n$ is either open or  
satisfies
$\mu^n(E)<\infty$ (or both) we have
\beq
\mu^n(E) = \sup\set{\mu^n(K): K\subset E,~K~\mbox{compact}}.
\eeq
\eit
These conditions are rather weak (for example, they are weaker than
inner- and outer-regularity), and are satisfied by most interesting
measures, including the counting measure on countable sets and the
Lebesgue measure on $\RR^n$ (see~\cite{rudin}, Theorem 2.14).

We say that a real-valued function $f$
defined on a metric space $(\calX,\rho)$
is {\em locally Lipschitz} if for each $x\in\calX$
there is an open $x\in U\subset\calX$
and a $0<C(x)<\infty$
such that 
\beq
\sup_{y\in U\setminus\set{x}}
\frac{\abs{f(x)-f(y)}}{\rho(x,y)}
&\leq& C(x).
\eeq
\hide{
We write
$\Lips^*(\calX,\rho)$ 
to denote
the set of all locally Lipschitz
on $f:\calX\to\RR$.
}
\hide{
In what follows,
if $E\subset\calX$ for some
metric space $(\calX,\rho)$,
we use the notation 
$E^c \defeq \calX\setminus E$
and
$\rho(x,E)=\inf\set{\rho(y,x):y\in E}$.
$C_c(\calX,\rho)$ will denote the space of the continuous functions
$f:\calX\to\RR$ with compact support.
}

\bethn
\label{thm:phinorm}
Let $\mu$ be a 
measure on a locally compact
metric
space 
$(\calX,\rho)$, where $\mu$
satisfies the regularity conditions (a)-(c) above.
Then for any $f\in L_1(\calX,\mu)$, $\phinorm{f}=0$ iff $f=0$ a.e. $[\mu]$.
\enthn
\bepf
Suppose $f\in L_1(\calX,\mu)$.
The case $f\leq 0$
a.e. $[\mu]$ 
is trivial, so we assume the existence of
a Borel $E\subset \calX$
such that
\beq
0<\mu(E)<\infty,
\qquad
f>0~\mbox{on}~E.
\eeq
Let $g(x)=\pred{x\in E}$ be the characteristic
function of $E$ and note that
$g\in L_1(\calX,\mu)$.

Theorems 2.24 and 3.14 in
\cite{rudin} (the first is Lusin's theorem)
provide a sequence of 
continuous functions
$h_n$
such that 
\beq
\sup_{x\in\calX} |h_n(x)| \leq \sup_{x\in\calX} |g(x)| = 1,
\qquad
\nrm{g-h_n}_1\to 0,
\eeq
which implies $h_n\to g$ a.e. $[\mu]$.
Thus by
Lebesgue's Dominated Convergence theorem, we have 
\beqn
\label{eq:fhn}
\iprod{f}{h_n}\to\iprod{f}{g}=\int_E f d\mu >0.
\eeqn

At this point we will need two facts:
\ben
\item continuous functions can be uniformly approximated by 
locally Lipschitz functions
\item locally Lipschitz functions can be uniformly approximated by
  finite linear combinations of members of
$\Lips(\calX,\rho)$
(defined in (\ref{eq:Lips});
\een
both are straightforward to verify.
It follows from (\ref{eq:fhn})
that 
the linear functional
$\iprod{f}{\cdot}$ cannot vanish on all of 
$\Lips(\calX,\rho)$, 
which implies 
$\phinorm{f}>0$.
\enpf

Theorem~\ref{thm:psinorm} shows that $\psinorm{\cdot}$ is
topologically equivalent to $\nrm{\cdot}_{L_1}$.
The norm strength of $\phinorm{\cdot}$ is a more interesting
matter.\hide{
We say that a metric $\rho$ is {\em
  sub-Hamming} if
\beq
\rho(x,y) &\leq& d_\ham(x,y)
\eeq
for $x,y\in\X^n$; note that conditions 
(\PSdiam)
and
(\PScons)
imply that every metric in a $\Psi$-dominated family is sub-Hamming.}
In the case of finite $\X$, 
$F_n=\ell_1(\X^n
)$ is a finite-dimensional space so all norms on
$F_n$ are trivially equivalent. 
Suppose $\X$ is a countable set 
(equipped with
the counting measure)
and $\rho$ has the property that
\beq
d = \inf_{x\neq y}\rho(x,y) > 0.
\eeq
The functions 
$g(x)=d\pred{f(x)>0}$
and
$h(x)=d\pred{f(x)<0}$
are both in
$\Lips(\X,\rho)$,
and since 
$d\nrm{f}_{1}=\abs{\iprod{f}{g}}+\abs{\iprod{f}{h}}$,
we have
\beqn
{\ts\oo2}d\nrm{f}_{1} \;\leq\;
\phinorm{f} \;\leq\;
\diam_{\rho}(\X) \nrm{f}_{1}
\eeqn
for all $f\in F_n$,
so the norms $\phinorm{\cdot}$ and $\nrm{\cdot}_{1}$ are equivalent
in this case.

Suppose, on the other hand, that
$T=\set{x_1, x_2, \ldots}$ forms a Cauchy sequence in 
the countable space
$\X$, with
$\delta_i = \rho(x_i,x_{i+1})$ approaching zero. Let
$f\in 
\ell_1(\X)
$ be such that
$f(x_{2k}) = -f(x_{2k-1})$ for $k=1,2,\ldots$
and $f(x)=0$ for $x\notin T$; then
\beqn
\label{eq:cauchysum}
\phinorm{f} \;\leq\; 
\sum_{k=1}^\infty \abs{f(x_{2k-1})}\delta_{2k-1}
\;\leq\;
\nrm{f}_{1}\sum_{k=1}^\infty \delta_{2k-1}.
\eeqn
If $\X=\QQ\cap[0,1]$ (the rationals in $[0,1]$)
 with $\rho(x,y)=|x-y|$ as the metric on $\X$, the
r.h.s. of (\ref{eq:cauchysum}) can be made arbitrarily small, so
for this metric space,
\beq
\inf\set{\phinorm{f} : \nrm{f}_{1}=1} &=& 0
\eeq
and $\phinorm{\cdot}$
is a strictly weaker norm than
$\nrm{\cdot}_{1}$.

Similarly, when $\X$ is a continuous set,
 $\phinorm{\cdot}$
will be strictly weaker than $\nrm{\cdot}_{L_1}$
in a fairly general setting.
As an example, take $n=1$, $\X=[0,1]$,
$\mu$ the Lebesgue measure on $[0,1]$,
and $\rho(x,y)=|x-y|$.
For $N\in\NN$, define
$\gamma_N:[0,1]\to\NN$ by
\beq
\gamma_N(x) &=& 
\max\set{0\leq k<N : k/N\leq x}.
\eeq
Consider the function
\beq
f_N(x) = (-1)^{\gamma_N(x)},
\eeq
for 
$N=2,4,6,\ldots$;
note that
$f$ is measurable and 
$\nrm{f}_{L_1}=1$.

For a fixed even $N$, 
define the $k$th segment
\beq
I_k = \set{x\in[0,1] : k\leq\gamma_N(x)\leq k+2}
= \sqprn{\frac{k}{N},\frac{k+2}{N}},
\eeq
for $k=0,2,\ldots,N-2$.
Since $\diam I_k = 2/N$, for any $g\in\Lips(\X,\rho)$, we have
\beq
\sup_{I_k} g(x) - \inf_{I_k} g(x) \leq 2/N;
\eeq
this implies
\beq
\int_{I_k} f_N(x) g(x) d\mu(x) \leq 2/N^2.
\eeq
Now $[0,1]$ is a union of $N/2$ such segments, so
\beq
\int_0^1 f_N(x) g(x) d\mu(x) \leq 1/N.
\eeq
This means that $\phinorm{f}$ can be made arbitrarily small
while $\nrm{f}_{L_1}=1$, so once again
and $\phinorm{\cdot}$
is a strictly weaker norm than
$\nrm{\cdot}_{L_1}$.

\end{appendix}

\section*{Acknowledgements}
I thank 
John Lafferty and
Kavita Ramanan for helpful discussions,
and Steven J. Miller for comments on the draft.

\end{document}